\documentclass[a4paper,11pt]{article}
\usepackage{latexsym}
\usepackage{amsfonts,amssymb}
\usepackage{amscd}
\usepackage[dvips]{graphicx}

\font\erm=cmr9

\author{A. K. Kwa\'sniewski}
\title{CAUCHY-RIEMANN EQUATIONS FOR CAYLEY NUMBERS' FUNCTIONS}

\newtheorem{defn}{Definition}

\newtheorem{lemma}{Lemma}

\begin{document}

\begin{center}
\noindent {\Large \textsc{Cauchy-Riemann Equations For Cayley Numbers' Functions}}  \\ 

\vspace{0.5cm}

\noindent A. Krzysztof Kwa\'sniewski 

\vspace{0.15cm}

\noindent {\erm the Dissident - relegated by Bia\l ystok University (*) authorities  }\\
\noindent {\erm from the Institute of Informatics to Faculty of Physics}\\
\noindent {\erm ul. Lipowa 41,  15 424  Bia\l ystok, Poland}\\

\noindent {\erm (*) former Warsaw University Division}\\
\noindent {\erm e-mail: kwandr@gmail.com}

\end{center}
% % % % % % % Abstract
\begin{abstract}
Since the discovery of octonions in 1843   by John T. Graves  [1] we seem to be still lacking a satisfactory if any  theory of  octave valued functions - satisfactory according to standard requirements or expectation from the side of  a theory like a one might look for.   Here is a proposal coming back to my twentieth century presentation of a perhaps nonstandard idea  hoping to be  coping with nonassociativity by an invention.

\end{abstract}

\vspace{0.4cm}


\vspace{0.1cm}

\noindent AMS Classification Numbers: 11S80,11R52,17D05   
\vspace{0.1cm}

\vspace{1cm}

By 1828 George Green  (born on 14 July 1793) had written his first and most important paper entitled  "An essay on the Application of Mathematical Analysis to the Theories of Electricity and Magnetism". 

In this essay , which runs to nearly 80 pages , George Green had formulated what today is called - Greens` Theorem. This theorem is used to derive immediately the Cauchy Theorem 

$$
	\oint_{\Gamma} f(z)dz = 0,
$$

\noindent as it leads to Cauchy-Riemann equations involved there.

So \emph{we now may celebrate the 180-th anniversary of this great achievement}.

In this note - (presented at "The Polish-Mexican Seminar on Generalized Cauchy -Riemann Structures and  Surface Properties of Crystals" - Kazimierz  Dolny; August 98 ) - 
 - one proposes the extension of the analyticity notion so that it includes also octonions and in general all composition algebras [2]. 

We also indicate the possibility of introducing the notion of analyticity for other
algebras (suggested by A. Z. Jadczyk in private communication).

\emph{The major aim of the note is to formulate the analyticity notion for an octonion algebra in a manner which would enable one to reestablish  the main theorems already known for quaternions and for Clifford  algebra valued regular functions.}

\vspace{0.6cm}

\textbf{I.} A composition algebra A is not necessarily associative, however, if it is, then
by means of the Cayley-Dickson procedure on can construct a new composition algebra
 $(A, \alpha)$ which is a direct sum of vector spaces $A\otimes A$ with the usual addition and multiplication by a real number while the internal product is defined by 

$$
	(x_1, y_1)(x_2, y_2) = (x_1 x_2 + \alpha y_2 y_1, x_1 y_2 + y_1 x_2), \ \ \ \  \alpha\neq 0,  a\in \mathbf{R}
$$

\noindent with the standard notation for conjugate elements in $A$. Conjugation in $(A, a)$ is defined by

$$
	\overline{(x, y)} = (\overline{x}, -y).
$$

In the following we restrict our discussion of analyticity concept to the more familiar case of ordinary composition algebras, i.e. complex numbers, quaternions and octonions which we shall call briefly just composition algebras, this being justified by the fact that most of our consideration are valid for the general case of any composition algebra .Let us start with a unified formulation of the algebras of complex numbers $\mathbf{C}$,  quaternions $\mathbf{C}$ and octonions $\Theta$.

From now on Greek indices $\mu, \nu, y, \sigma, ...$ will run  from $0$ to $1$ for $\mathbf{C}$, from $0$ to $3$ for $\mathbf{Q}$ and from $0$ to $7$ for $\mathbf{\Theta}$, while Latin indices $i, j, k, ...$ shall take correspondingly values $1$, or $l, 2, 3$, or $l, 2, ..., 7$ (summation convention is used). The algebras $\mathbf{C}$, $\mathbf{Q}$, $\mathbf{\Theta}$ can be defined via

\begin{equation}\label{eq:1I}
	e_\mu e_\nu = c^{\sigma}_{\mu \nu} e_\sigma
\end{equation}

\noindent where

\begin{equation}\label{eq:2I}
	c^{\sigma}_{ij} = -\delta_{ij}\delta^\sigma_0 + \epsilon_{ij}\delta^\sigma_k, \ \ \ \ 
	c^\sigma_{0\mu} = c^\sigma_{\mu 0} = \delta^\sigma_\mu 
\end{equation}

\begin{figure}[ht]
\begin{center}
	\includegraphics[width=60mm]{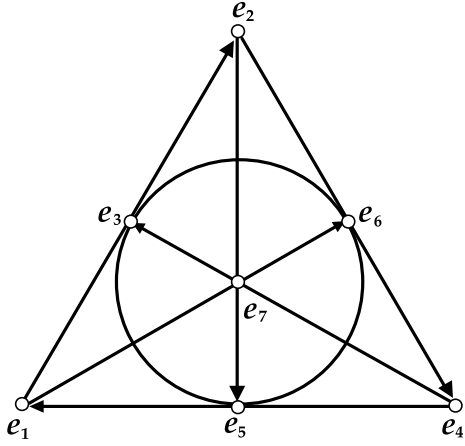}
	\caption{$7$-th element - plane, projective geometry  picture \cite{4} \label{fig:1}}
\end{center}
\end{figure}

\vspace{0.2cm}
\noindent For $i = j = k = 1$ we have $\mathbf{C}$, for $i, j, k = 1, 2, 3$ we get $\mathbf{Q}$ and for $i, j, k = 1, . . . , 7$    we shall obtain the algebra of octonions (see \textbf{Fig.\ref{fig:1}}) if one defines for which triples $(i, j, k)$  $\epsilon_{ijk}= + 1$. In the case of octonions we must add to (\ref{eq:2I}) the following conditions

\begin{equation}\label{eq:3I}
	\epsilon_{123} = \epsilon_{145} = \epsilon_{176} = \epsilon_{246} = \epsilon_{347} = \epsilon_{536} = \epsilon_{725} = + 1.
\end{equation}

\begin{center}
The rules: \\
$e_1 e_3 = e_2, e_2 e_6 = e_4, e_4 e_5 = e_1, e_3 e_6 = e_5, e_1 e_7 = e_6, e_2 e_7 = e_5, e_4 e_7 = e_3$.
\end{center}

Let now A be any of the algebras $\mathbf{C}$, $\mathbf{Q}$ or $\mathbf{\Theta}$, then $x\in A$ can be represented as

$$
	x = x^\mu e_\mu, \ \ \ x_\mu \in \mathbf{R},
$$
\noindent where
\begin{equation}\label{eq:4I}
	x^\mu = \frac{1}{2}\left( x \overline{e}_\mu, + e_\mu \overline{x} \right)
\end{equation}
\noindent and
$$
	e_0 = \overline{e}_0, \ \ \ \ e_i = -\overline{e}_i
$$

\noindent The trace, a linear mapping of $A$ into $\mathbf{R}$, is defined then by

\begin{equation}\label{eq:5I}
	\mathrm{Tr} \ e_i = 0, \ \ \ \ \mathrm{Tr} \ e_0 = N,
\end{equation}

\noindent where $N = \mathrm{dim} \ A$. Using this trace mapping one may introduce a scalar product in $A$

\begin{equation}\label{eq:6I}
	x, y \in A, \ \ \ \ \langle x | y \rangle = \frac{1}{N} \mathrm{Tr}\left( \overline{x} y \right),
\end{equation}

\noindent which has the property

\begin{equation}\label{eq:7I}
	\langle xy | xy \rangle = \langle x | x \rangle \langle y | y \rangle
\end{equation}

\noindent Using the definition (\ref{eq:2I}) of the $A$ algebra structure constants one may derive
the following properties of $c^\sigma_{\mu \nu}$:

\begin{equation}\label{eq:8I}
	c^0_{\mu \nu} = g_{\rho \nu} \ \ \ \ \mathrm{where} \ \ \ \ g_{\rho \nu} = \left\{ 
		\begin{array}{rl}
			\delta_{0\nu}, & \rho = 0 \\
			-\delta_{i\nu}, & \rho = i
		\end{array}
	\right.
\end{equation}

\begin{equation}\label{eq:9I}
	c^\sigma_{\sigma 0} = N, \ \ \ \ c^\sigma_{\sigma k} = 0,
\end{equation}

\begin{equation}\label{eq:10I}
	c^k_{ij} = \epsilon_{ijk},
\end{equation}

\begin{equation}\label{eq:11I}
	c^\sigma_{\mu \nu} g_{\sigma \rho} = c^\sigma_{\rho \mu} g_{\sigma \nu},
\end{equation}

\noindent or

$$
	c^{\overline{\rho}}_{\mu \nu} = c^{\overline{\nu}}_{\rho \nu} \ \ \ \
	\mathrm{where} \ \ \ \ c^{\overline{\rho}}_{\mu \nu} = \left\{ 
		\begin{array}{rl}
			c^0_{\mu\nu}, & \rho = 0 \\
			-c^i_{\mu\nu}, & \rho = i
		\end{array}
	\right.
$$

\begin{equation}\label{eq:12I}
	e_\mu \overline{e}_\nu  +  e_\nu \overline{e}_\mu  = 2 \delta_{\mu\nu} e_0.
\end{equation}

With the help of (\ref{eq:8I}-\ref{eq:11I}) one can prove an important lemma. For that to do 
let us introduce a differential linear operator of the form 

$$
	L = L^\sigma e_\sigma
$$

\noindent acting on $A$ and a mapping 

$$
	U : A \rightarrow A
$$

\noindent for which $U^\mu(x)$ ($U(x) = U^\mu(x) e_\mu$) are real differentiable functions \cite{2}.

\begin{lemma}
\begin{equation}\label{eq:13I}
	\mathrm{Tr} \left\{ L(U_q) \right\} = \mathrm{Tr} \left\{ (LU)_q \right\}, \ \ \ \ \forall q \in A.
\end{equation}
\end{lemma}

\noindent Another useful lemma can be established using the relations (\ref{eq:8I}-\ref{eq:10I}) \cite{2}.

\begin{lemma}
\begin{equation}\label{eq:14I}
	c^\nu_{\sigma\mu} \partial^\sigma \equiv 2\delta^\mu_0 \partial_\nu - c^\sigma_{\nu\mu} \partial_\sigma,
\end{equation}
\noindent where $\partial^\sigma \equiv \partial / \partial x_\sigma$.
\end{lemma}

\vspace{0.6cm}
\textbf{II.} 
In this section we construct a matrix representation of $A = \mathbf{C}, \mathbf{Q}, \mathbf{\Theta}$ with
usual addition and multiplication of matrices as operations in A. The matrices
will have \textbf{operator entries} as one of specifications of $A$ is nonassociative. In
the associative cases the operator entries simply become matrices.
It is well known that $\mathbf{C}$ can be isomorphically represented by a set of matrices
of the form

\begin{equation}\label{eq:1II}
	\mathbf{C} \in z, \ \ \ \ z = \Bigg[
	\begin{array}{rr}
		x & -y \\
		y &  x
	\end{array}
	\Bigg] , \ \ \ \ 
	x, y \in \mathbf{R}
\end{equation}

\noindent The conjugation $\sigma:\mathbf{C} \rightarrow \mathbf{C}$ can also be represented via matrix multiplication in the following way:

\begin{equation}\label{eq:2II}
	\sigma(z) = jzj, \ \ \ \ \mathrm{where} \ \ \ \ j^2 = id, \ \ \ \, jz = \overline{z} j
\end{equation}

\noindent due to commutativity of $\mathbf{C}$. The $j$ matrix is of the form

\begin{equation}\label{eq:3II}
	j = \Bigg[
	\begin{array}{rr}
		1 & 0 \\
		0 & -1
	\end{array}
	\Bigg] ,
\end{equation}

To work out a similar construction for quaternions it is sufficient to notice the
essence of the above representation which was the Cayley-Dickson procedure  applied to $\mathbf{R}$.
$\mathbf{Q}$ can be then isomorphically represented as a set of matrices of the form

\begin{equation}\label{eq:4II}
	q \in \mathbf{Q}, \ \ \ \ q = \Bigg[
	\begin{array}{cc}
		z_1 & -z_2 j \\
		z_2 j & z_1
	\end{array}
	\Bigg], \ \ \ \
	z_1, z_2 \in \mathbf{C}
\end{equation}

\noindent This form again is a manifestation of the Cayley-Dickson procedure.  

\noindent The conjugation $\sigma(q) = \overline{q}$ also can be realized by

$$
	\sigma(q) = \epsilon \circ q \circ \epsilon, \ \ \ \
	e^2 = \mathrm{id}, \epsilon \cdot q = \overline{q} \cdot \epsilon
$$

\noindent but this time, because of noncommutativity of $\textbf{Q}$, $\epsilon \circ q \circ \epsilon$ does not denote simply
matrix multiplication as we must have

$$
	\epsilon \circ (q_1 \circ q_2) \circ \epsilon = (\epsilon \circ q_2 \circ \epsilon)(\epsilon \circ q_1 \circ \epsilon).
$$

However, it is enough to say that the \emph{ \textbf{$\epsilon$ operator acts by matrix multiplication} under
the condition  that  this  multiplication reverses  the order  of the  product  of
$q$-matrices whenever they are sandwiched between two $\epsilon$ operators} (matrices).
With this in mind

\begin{equation}\label{eq:5II}
	\epsilon = \Bigg[
	\begin{array}{rr}
		j & 0 \\
		0 & -j
	\end{array}
	\Bigg].
\end{equation}

\noindent Let us identify

$$
	z_1 \in \mathbf{Q}, \ \ \ \
	z_1 = \Bigg[
	\begin{array}{rr}
		z_1 & 0 \\
		0 & z_1
	\end{array}
	\Bigg], \ \ \ \
	z_1 \in \mathbf{C}.
$$

\noindent Then

\begin{equation}\label{eq:6II}
	q \in \mathbf{Q}, \ \ \ \
	q = z_1 + z_2 j,
\end{equation}

\noindent where

$$
	j = \Bigg[
	\begin{array}{rr}
		0 & -j \\
		-j & 0
	\end{array}
	\Bigg] \in \mathbf{Q}.
$$

\noindent Now the action of $\epsilon$ on $\mathbf{Q}$ can be defined as follows:

\begin{equation}\label{eq:7II}
	\epsilon \circ q \circ \epsilon  =  \epsilon z_1 \epsilon + \epsilon j z_2 \epsilon \equiv \overline{q}
\end{equation}

\noindent and as

$$
	\epsilon j = -j \epsilon, \ \ \ \
	z j = j \overline{z}
$$

\noindent one has

$$
	\epsilon \circ q = \overline{q} \circ \epsilon.
$$

\noindent The natural representation of imaginary units $e_1, e_2, e_3 \in \mathbf{Q}$ is then given by

\begin{equation}\label{eq:8II}
	e_1 = \Bigg[
	\begin{array}{rr}
		i & 0 \\
		0 & i
	\end{array}
	\Bigg], \ \ \ \
	e_2 = \Bigg[
	\begin{array}{rr}
		0 & -j \\
		-j & 0
	\end{array}
	\Bigg], \ \ \ \
	e_3 = \Bigg[
	\begin{array}{rr}
		0 & -ij \\
		ij & 0
	\end{array}
	\Bigg]
\end{equation}

Note that the product of  $q_1, q_2 \in \mathbf{Q}$ is realized by the usual matrix multiplication
$(q_1 \circ q_2 = q_1 q_2)$ while $\epsilon \notin \mathbf{Q}$.  

$\epsilon$ is a specific operator acting on $\mathbf{Q}$, $\epsilon$ is in a sense a "square root" of the conjugation operator s and can be thought of as the matrix (\ref{eq:5II})  but then one must remember that though it acts by matrix multiplication ...   it  reverses the order of $q$-matrices  -  if sandwiched between two $\epsilon$  matrices.
Similarly to previous cases, octonions can be represented by

\begin{equation}\label{eq:9II}
	\theta \in \Theta, \ \ \ \
	\theta = \Bigg[
	\begin{array}{lc}
		q_1 & -q_2 \epsilon \\
		q_2 \epsilon & q_1
	\end{array}
	\Bigg] \equiv q_1 + q_1 E,
\end{equation}

\noindent where

$$
	E = \Bigg[
	\begin{array}{rr}
		0 & -\epsilon \\
		\epsilon & 0
	\end{array}
	\Bigg], \ \ \ \
	q_1, q_2, E \in \Theta.
$$

\noindent The representation is given once the multiplication law in $\Theta$ is defined. It is given by

\begin{equation}\label{eq:10II}
	\begin{array}{ll}
	(a) &  (q_1 \epsilon) \circ (q_2 \epsilon) = \epsilon \circ (\overline{q}_1 q_2) \circ \epsilon = \overline{q}_2 q_1, \\
	(b) &  (q_1 \epsilon) \circ q_2 = (q_1 \overline{q}_2) \epsilon, \\
	(c) &  q_1 \circ (q_2 \epsilon) = (q_2 q_1) \epsilon \\
	(d) & \epsilon \cdot \epsilon = \epsilon \epsilon = 1.
	\end{array}
\end{equation}

The rules (\ref{eq:10II}) can be derived (for that form of matrix representation (\ref{eq:9II})) from the Moufang identities \cite{3}. For example (\ref{eq:10II}a) can be derived from

$$
	x(yz)x = (xy)(zx) \ \ \ \ 
	\forall x,y,z \in A; \ \ \ \
	\mathrm{where} \ A  \mathrm{-alternative \ algebra}.
$$

\vspace{0.2cm}
To end up: multiplication in $\mathbf{\Theta}$ is just matrix multiplication where the rules of dealing with expressions involving $\epsilon$ symbols  are given by (\ref{eq:10II}). These rules apply to any specification of $A$. The "square root" $\epsilon$ of the conjugation operator $\sigma$ (with respect to (\ref{eq:10II}) multiplication) is given by the matrix

\begin{equation}\label{eq:11II}
	\epsilon = \Bigg[
	\begin{array}{rr}
		\epsilon & 0 \\
		0 & -\epsilon
	\end{array}
	\Bigg]
\end{equation}

\noindent and again $\epsilon \notin \mathbf{\Theta}$; it is an operator acting on $\mathbf{\Theta}$ similarly as $E$ does on $\mathbf{Q}$. In the natural representation, generators of $\mathbf{\Theta}$ have the form

\begin{equation}
	{
	e_1 = \Bigg[
	\begin{array}{ll}
		e_1 & 0 \\
		0 & e_1
	\end{array}
	\Bigg], \ \ \
	e_2 = \Bigg[
	\begin{array}{ll}
		e_2 & 0 \\
		0 & e_2
	\end{array}
	\Bigg], \ \ \
	e_3 = \Bigg[
	\begin{array}{ll}
		e_3 & 0 \\
		0 & e_3
	\end{array}
	\Bigg], \ \ \
	e_4 = \Bigg[
	\begin{array}{lc}
		0 & -\epsilon \\
		\epsilon & 0
	\end{array}
	\Bigg], \ \ \
	\atop
	e_5 = \Bigg[
	\begin{array}{lc}
		0 & -e_1\epsilon \\
		e_1\epsilon & 0
	\end{array}
	\Bigg], \ \ \
	e_6 = \Bigg[
	\begin{array}{lc}
		0 & -e_2\epsilon \\
		e_2\epsilon & 0
	\end{array}
	\Bigg], \ \ \
	e_7 = \Bigg[
	\begin{array}{lc}
		0 & -e_3\epsilon \\
		e_3\epsilon & 0
	\end{array}
	\Bigg].
	}
\end{equation}

\noindent Using this representation one proves the following \cite{2}

\begin{lemma}
$$
	\forall \ \Theta, u \in A \ \ \ \ \
	\Theta (\overline{\Theta}u) = (\Theta \overline{\Theta}) u.
$$
\end{lemma}

\emph{Introducing then the operators defined on functions on $Q$}
$$
	\partial_{\overline{q}} \equiv e^\mu \partial_\mu
$$

\noindent where $\mathbf{Q} \ni q = x^\mu e_\mu$ and $\partial_q \equiv \overline{e}^\mu \partial_\mu$ \emph{we have for octonions}

\begin{equation} \label{eq:13II}
	\partial_\theta = \Bigg[
	\begin{array}{cc}
		\partial_{q_1} & \partial_{\overline{q}_2} \epsilon \\
		-\partial_{\overline{q}_2} \epsilon & \partial_{q_1}
	\end{array}
	\Bigg], \ \ \ \ 
	\partial_{\overline{\theta}} = \Bigg[
	\begin{array}{cc}
		\partial_{\overline{q}_1} & \partial_{{q}_2} \epsilon \\
		-\partial_{{q}_2} \epsilon & \partial_{\overline{q}_1}
	\end{array}
	\Bigg], \ \ \ \ 
\end{equation}

\begin{equation}\label{eq:14II}
	\partial_{\overline{\theta}} \circ \partial_{{\theta}} = \partial_{{\theta}} \circ \partial_{\overline{\theta}} = 
	\Bigg[
	\begin{array}{cc}
		\diamondsuit_\infty & 0 \\
		0 & \diamondsuit_\infty
	\end{array}
	\Bigg],
\end{equation}

\noindent where $\Box_8 \equiv \diamondsuit_\infty = \partial^\mu \partial_\mu$ and the $\infty$-sign stands for   horizontal $8$, because of my editorial limitations.

\vspace{0.6cm}
\textbf{III.} The Cauchy-Riemann (C-R) equations for $\mathbf{C}$ can be written in the form

$$
	\Bigg[
	\begin{array}{rr}
		\partial_0 & -\partial_1 \\
		\partial_1 & \partial_0
	\end{array}
	\Bigg]
	\Bigg[
	\begin{array}{rr}
		u_0 & -u_1 \\
		u_1 & u_0
	\end{array}
	\Bigg] = 0
$$

\noindent or $\partial_{\overline{z}}U = 0$, where $\partial_\mu = \partial/\partial x^\mu$, $U = U^\mu e_\mu$; $\mu= 0,1$.

\vspace{0.4cm}
This definition of analytic function $U$ can be extended to any algebra $A = \mathbf{C}, \mathbf{Q}, \mathbf{\Theta}$.
Let $U$ be an $A$-valued function on $A$ with $U^\mu (x_0, x_1, ..., x_{N-1})$ functions differentiable with respect to $x_\nu$. The $U$ can be represented as

$$
	U = \Bigg[
	\begin{array}{cc}
		a & -b\alpha \\
		b\alpha & a
	\end{array}
	\Bigg]
$$

\noindent with $\alpha = 1, j, \epsilon$ correspondingly to the chosen case; $A = \mathbf{C}, \mathbf{Q}, \mathbf{\Theta}$.

\begin{defn}
U is called left \textbf{$A$-analytic} iff
\end{defn}

\begin{equation}\label{eq:1III}
	\partial_{\overline{\theta}}U = 0 \ \ \ \ \mathrm{or} \ \ \ \ 
	\Bigg[
	\begin{array}{cc}
		\partial_{\overline{q}_1} & -\partial_{{q}_2}\epsilon \\
		\partial_{{q}_2}\epsilon & \partial_{\overline{q}_1}
	\end{array}
	\Bigg]
	\Bigg[
	\begin{array}{cc}
		a & -b\alpha \\
		b\alpha & a
	\end{array}
	\Bigg] = 0
\end{equation}

\begin{defn}
$U$ is called left \textbf{$A$-anti}analytic iff
\end{defn}

\begin{equation}\label{eq:2III}
	\partial_{{\theta}}U = 0 \ \ \ \ \mathrm{or} \ \ \ \ 
	\Bigg[
	\begin{array}{cc}
		\partial_{{q}_1} & \partial_{\overline{q}_2}\epsilon \\
		-\partial_{\overline{q}_2}\epsilon & \partial_{{q}_1}
	\end{array}
	\Bigg]
	\Bigg[
	\begin{array}{cc}
		a & -b\alpha \\
		b\alpha & a
	\end{array}
	\Bigg] = 0
\end{equation}

\noindent As a conclusion from Section II we get  \cite{2}

\begin{lemma}
An $A$-analytic or $A$-antianalytic function is a harmonic function, i.e.
$$
	\diamondsuit_N U = 0.
$$
\end{lemma}

There exists a lot of \emph{$A$-analytic} functions. The infinite number of examples is given by simple combinations of \emph{$B$-analytic} and \emph{$B$-antianalytic} functions where  $A = (B, - 1)$ (see Section I).
The Cayley-Dickson procedure  inherent in this definition allows us to relate octonion analy-ticity to quaternion or via quaternion to complex, "usual" analyticity. Let us consider in more detail octonion-antianalyticity as an example. The octonion function can be written in three equivalent forms

$$
	\begin{array}{c}
	U = a + bE, \ \ \ \ \ a, b, E \in \Theta, \\
	U = A_1 + B_1 j + (A_2 + B_2 j) E, \ \ \ \ A_1, A_2, B_1, B_2 \in \mathbf{C}, \\
	U = U^\mu e_\mu, \ \ \ \ U^0, U^1, ..., U^\mu \in \mathbf{R}.
	\end{array}
$$

\noindent Introducing the notation

$$
\begin{array}{cc}
	\partial_\theta = \Bigg[
	\begin{array}{cc}
		\partial_{q_1} & \partial_{\overline{q}_2}\epsilon \\
		-\partial_{\overline{q}_2}\epsilon & \partial_{q_1}
	\end{array}
	\Bigg], & \partial_{q_1} = \Bigg[
	\begin{array}{cc}
		\partial_{y_1} & \partial_{\overline{z}_1}j \\
		-\partial_{\overline{z}_1}\epsilon & \partial_{y_1}
	\end{array}
	\Bigg]
	\\
	\partial_{y_1} = \Bigg[
	\begin{array}{cc}
		\partial_0 & \partial_1 \\
		-\partial_1 & \partial_0
	\end{array}
	\Bigg],
	&
	\partial_{\overline{z}_1} = \Bigg[
	\begin{array}{cc}
		\partial_2 & -\partial_3 \\
		\partial_3 & \partial_2
	\end{array}
	\Bigg], ..., 
	\partial_{z_2} = \Bigg[
	\begin{array}{cc}
		\partial_6 & \partial_7 \\
		-\partial_7 & \partial_6
	\end{array}
	\Bigg] ...
\end{array}
$$

\noindent we can write \textbf{octonionic antianalyticity C-R euqations} in three equivalent forms:

\vspace{0.4cm}
\noindent \textbf{Quaternionic form:}

\begin{equation}\label{eq:3III}
{
	\partial_{q_1}a + \overline{b}\partial^l_{\overline{q}_2} = 0
\atop
	-\partial_{\overline{q}_2} \overline{a} + b \partial^l_{q_1} = 0,
}
\end{equation}

\noindent where $\partial^l$ means action to the left.

\vspace{0.4cm}
\noindent \textbf{Complex form:}

\begin{equation}\label{eq:4III}
	\begin{array}{r}
		\partial_{y_1}A_1  +  \partial_{\overline{z}_1}B  +  \partial_{\overline{y}_2}\overline{A}_2  +  \partial_{z_2}B_2 = 0 \\
		\partial_{y_1}B_1  -  \partial_{\overline{z}_1}A_1  -  \partial_{{y}_2}{B}_2  +  \partial_{\overline{z}_2}\overline{A}_2 = 0 \\
		-\partial_{\overline{y}_2}A_1  -  \partial_{\overline{z}_2}\overline{B}_1  +  \partial_{{y}_1}{A}_2  +  \partial_{{z}_1}{B}_2 = 0 \\
		-\partial_{\overline{z}_2}A_1  +  \partial_{\overline{y}_2}{B}_1  +  \partial_{\overline{y}_1}{B}_2  -  \partial_{\overline{z}_1}{A}_2 = 0 \\
	\end{array}
\end{equation}

\noindent \textbf{Real form:}

\begin{equation}\label{eq:5III}
\left[
\begin{array}{rrrrrrrr}
	\partial_0 & \partial_1 & \partial_2 & \partial_3 & \partial_4 & \partial_5 & \partial_6 & \partial_7 \\
	-\partial_1 & \partial_0 & \partial_3 & -\partial_2 & \partial_5 & -\partial_4 & -\partial_7 & \partial_6 \\
	-\partial_2 & -\partial_3 & \partial_0 & \partial_1 & \partial_6 & \partial_7 & -\partial_4 & -\partial_5 \\
	-\partial_3 & \partial_2 & -\partial_1 & \partial_0 & \partial_7 & -\partial_6 & \partial_5 & -\partial_4 \\
	-\partial_4 & -\partial_5 & -\partial_6 & -\partial_7 & \partial_0 & \partial_1 & \partial_2 & \partial_3 \\
	-\partial_5 & \partial_4 & -\partial_7 & \partial_6 & -\partial_1 & \partial_0 & -\partial_3 & \partial_2 \\
	-\partial_6 & \partial_7 & \partial_4 & -\partial_5 & -\partial_2 & \partial_3 & \partial_0 & -\partial_1 \\
	-\partial_7 & -\partial_6 & \partial_5 & \partial_4 & -\partial_3 & -\partial_2 & \partial_1 & \partial_0 \\
\end{array}
\right]
\left[
\begin{array}{l}
U_0 \\
U_1 \\
U_2 \\
U_3 \\
U_4 \\
U_5 \\
U_6 \\
U_7
\end{array}
\right] = \vec{O}
\end{equation}

Left octonion-analyticity C-R eqs. are obtained by replacing $a$; by $-a$; in (\ref{eq:3III}).
Correspondingly, complex and real forms of quaternion-analyticity conditions  in the analogous notation 
($u = u^\mu e_\mu$)

$$
	U = A + Bj, \ \ \ \
	y = x_0 + x_1 i, \ \ \ \
	q = y + zj \ \ \ \
	z = x_2 + x_3 i,
$$

\noindent are given by:

\begin{equation}\label{eq:6III}
	\partial_{\overline{y}} A  -  \partial_{\overline{z}} \overline{B}  =  0, \ \ \ \ 
	\partial_{\overline{z}} \overline{A}  -  \partial_{\overline{y}} B  =  0, \ \ \ \ 
	\mathrm{(complex \ form)},
\end{equation}

\begin{equation}\label{eq:7III}
\left[
\begin{array}{rrrr}
	\partial_0 & -\partial_1 & -\partial_2 & -\partial_3 \\
	\partial_1 &  \partial_0 & -\partial_3 &  \partial_2 \\
	\partial_2 &  \partial_3 &  \partial_0 & -\partial_1 \\
	\partial_3 &  \partial_2 &  \partial_1 &  \partial_0
\end{array}
\right]
\left[
\begin{array}{l}
	U_0 \\
	U_1 \\
	U_2 \\
	U_3 
\end{array}
\right] = 0 \ \ \ \
	\mathrm{(real \ form).}
\end{equation}

\vspace{0.4cm}
The above formulation of analyticity coincides for $A = \mathbf{C}$ with Cauchy-Riemann and for $A = \mathbf{Q}$ with Fueter's analyticity.  \\
To see the latter we shall write C-R  equations for quaternions in another form.

\vspace{0.2cm}
\noindent Let us introduce the notation

$$
	\vec{U} = (U_1, U_2, U_3), \ \ \ \ \vec{V} = (\partial_1, \partial_2, \partial_3)
$$

\noindent for quaternions and

$$
	\vec{U} = (U_1, U_2, ..., U_7), \ \ \ \ \vec{V} = (\partial_1, \partial_2, ..., \partial_7)
$$

\noindent for octonions. Then (\ref{eq:7III}) can be written in the form

\begin{equation}\label{eq:8III}
	\partial_0 U_0 = \vec{V} \cdot \vec{U}, \ \ \ \ \partial_0\vec{U} = -\vec{V}U_0 - \vec{V} \times \vec{U},
\end{equation}

\noindent while (\ref{eq:5III}) is equivalent to

\begin{equation}\label{eq:9III}
	\partial_0 U_0 = \vec{V} \cdot \vec{U}, \ \ \ \ \partial\vec{U} = -\vec{V}U_0 - \vec{V} \otimes \vec{U},
\end{equation}

\noindent where the "\textbf{octonionic vector product $\otimes$}" is defined by

\begin{equation}\label{eq:10III}
	(\vec{V} \otimes \vec{U})_j = \sum_{i,k} \epsilon_{jki} \partial_k U_i
\end{equation}

\noindent with $\epsilon_{jki}$ satisfying (\ref{eq:3I}).

One component of the $\otimes$-vector product is an algebraic sum of six terms because the $(k, i)$ pair index takes six values for an index $j$ being fixed.
A more straightforward real form of C-R equations for $\mathbf{C}$, $\mathbf{Q}$ or $\mathbf{\Theta}$ is the equation

\begin{equation}\label{eq:11III}
	C^\rho_{\nu\sigma} \partial^\nu U^\sigma = 0
\end{equation}

\noindent equivalent to (\ref{eq:1III}). This however does not exhibit the structure originating from the Cayley-Dickson procedure. In the representation of $e_\mu$, we have given before, $\partial_{\overline{\Theta}} = e^\mu \partial_\mu$ acting as a linear operator on $A$, can be represented by a matrix in $a\{e_\mu\}$ basis. In view of Lemma 2 of Section I the matrix elements of this operator are given by the expression

$$
	\left( \partial_{\overline{\Theta}} \right)^\mu_\nu = C^\mu_{\sigma\nu} \partial_\mu  -  C^\sigma_{\nu\mu} \partial_\sigma .
$$

\vspace{0.6cm}
\textbf{IV.}  In this section we introduce the definition of analyticity for any algebra with unit element as was proposed by A. Z. Jadczyk (private communication); then we show that for an ordinary composition algebras it is exactly the same notion as the one we have introduced in previous sections. Let $A$ be now any algebra with unit element and let $U$ be a differentiable mapping $U: A \rightarrow A$. The derivative of $U$ at $x \in A$ is then an $R$-linear transformation $U'_x$ of $A$ and can be written

\begin{equation}\label{eq:1IV}
	U'_x(h) = h^\sigma \partial_\sigma U^\mu e_\mu ,
\end{equation}

\noindent where $A \ni h = h^\sigma e_\mu$, $U^\mu e_\mu = U \in A$, $U^\nu \in R$ and $\{e_\nu\}$ is a basis of $A$.

\vspace{0.4cm}
C-R equations are conditions on such a linear transformation $U'_x$; conditions related to the algebraic structure of $A$. The requirement of $A$-linearity though seemingly natural is very naive and yields a notion void of content already for quaternions [2].

A. Z. Jadczyk proposed a weakened condition. Let $f$ be any linear mapping $f: A \rightarrow A$. A trace of that mapping is then defined by

\begin{equation}\label{eq:2IV}
	\mathrm{Tr} f = f(e_\mu)^\mu ,
\end{equation}

\noindent where $\{e_\mu\}$ is a certain basis of $A$ and for  $a\in A$, $a^\mu$ denotes the $\mu$-th coordinate of a in a basis $\{e_\mu\}$.  Then as a generalization of C-R equations he proposes the equality of traces of the following two linear mappings

\begin{equation}\label{eq:3IV}
{
	f_{1,q}: A \rightarrow A, \ \ \ \ 
	f_{1,q}(a) = U'_x(qa),
\atop
	f_{2,1}: A \rightarrow A, \ \ \ \
	f_{2,q}(a) = U'_x(q)a, \ \ \ \
	q,a \in A.
}
\end{equation}

\vspace{0.2cm}
\noindent Generalized C-R eqs. for an algebra $A$ then have the form

\begin{equation}\label{eq:4IV}
	\forall q \in A \ \ \ \
	\left\{ U'_x(qe_\mu) \right\}^\mu = \left\{ U'_x(q)e_\mu \right\}^\mu 
\end{equation}

\noindent or equivalently:

\begin{equation}\label{eq:5IV}
	\left\{ U'_x(e_\mu e_\nu) \right\}^\nu = \left\{ U'_x(e_\mu)e_\nu \right\}^\nu .
\end{equation}

\noindent One then easily finds the C-R eqs. in terms of structure constants

\begin{equation}\label{eq:6IV}
	C^\sigma_{\mu \nu} \partial_\sigma U^\nu  =  C^\sigma_{\nu \sigma} \partial_\mu U^\nu
\end{equation}

\noindent In what follows we show that this definition for $A = \mathbf{C}, \mathbf{Q}, \mathbf{\Theta}$ is equivalent to  (\ref{eq:11III}).

\vspace{0.4cm}
At first let us notice that the definition (\ref{eq:2IV}) of the trace of a linear mapping for composition algebras coincides with that of Section I because of (\ref{eq:9I}). Secondly, (\ref{eq:6IV}) educes for composition algebras to

\begin{equation}\label{eq:7IV}
	C^\sigma_{\mu \nu} \partial_\sigma U^\nu  =  N \partial_\mu U^0
\end{equation}

\noindent and this definition is equivalent to ours because of Lemma 2 of Section l. The fact that we have a factor $N$ instead of 2 on the right-hand side of (\ref{eq:7IV}) is not important in view of the lemma (due to A. Z. Jadczyk).

\begin{lemma}
Let $F_\aleph$ denote the vector space of functions $f : A \rightarrow A$ satisfying the $\kappa-C-R$ eqs.
\begin{equation}\label{eq:8IV}
	C^\sigma_{\mu \nu} \partial_\sigma f^\nu  =  \kappa \partial_\mu f^0,
\end{equation}
where $\kappa, \kappa' \in \mathbf{R}$ and  $A$  is a composition algebra. Then for $\kappa, \kappa' \in \mathbf{R}, \kappa\kappa' - 1 \neq 0$ there
exists an isomorphism  $T_{\kappa\kappa'} : F_\kappa \rightarrow F_{\kappa}$. 
\end{lemma}

\noindent \emph{Proof:} The isomorphism is defined by $(T_{\kappa\kappa'}f)^i = f^i, (T_{\kappa\kappa'}f)^0 = \mu f^0$,
where $\mu = (1 - \kappa)/(1 - \kappa\kappa')$. One may check that $T_{\kappa\kappa'}f$ satisfies (\ref{eq:8IV}) with $\kappa'$ instead of $\kappa$  (use \ref{eq:9I}). 

\vspace{0.2cm}
\emph{Remark.} For $\kappa = N$ as in (\ref{eq:7IV}) one should take $\mu = (N -1)/(2N -1)$ to get $\kappa' = 2$ in (\ref{eq:8IV}) which then coincides with (\ref{eq:11III}) because of Lemma 2, Section \textbf{I}.

\vspace{0.4cm}
We have already argued in Sec. \textbf{III} that there are many examples of $A$-analytic functions, where $A$ is a composition algebra. However, this set of functions does not include $U(x) = x^n (n > 1)$ functions (except for $N = 2$) although it does include other $R$-homogeneous functions of degree $n$ and these $A$-analytic homogeneous functions play the role similar to $x^n$ in complex analytic functions theory .

To illustrate the above statement we quote \cite{1} the following

\begin{lemma}
Let $x \in A = \mathbf{C}, \mathbf{Q}, \mathbf{\Theta}; U(x) = x^2$; then $U$ satisfies C-R eqs. iff $N = 2$.
\end{lemma}

\vspace{0.4cm}
\noindent \textbf{To end up let us make three remarks.}

\vspace{0.4cm}
\textbf{1.} There exists a formulation of C-R eq. based on the analogy with Clifford  algebra product.   Let $U$ be a function on $A = Q, \Theta  (U^\mu e_\mu = U)$. We shall call $U_0$ the scalar part and $\vec{U} = (U_1, ..., U_{N-1})$ the vector part of $U$ and we shall represent $U$ by a pair $U  = (U_0, \vec{U})$ and similarly $\partial_{\overline{x}} = (\partial_0, \vec{\nabla})$, $x \in A$. \\
Then  $\kappa-C-R$ eqs. for $\kappa = 2$ can be written in the form

\begin{equation}\label{eq:9IV}
	\left( \partial_0, \vec{\nabla} \right)   \circ   \left( U_0, \vec{U} \right) = 0,
\end{equation}

\noindent where

$$
	(a_0,\vec{a})  \circ  (b_0,\vec{b})  =  (a_0 b_0 - \vec{a}\vec{b}, a_0\vec{b} + \vec{a}b_0 + \vec{a} \otimes \vec{b});
$$

\noindent \emph{$\otimes$ denotes the octonionic vector product} coinciding with the usual one for the quaternionic subalgebra of $\Theta$. One easily notices that eqs. (\ref{eq:9III}-\ref{eq:10III}) are just specifications of (\ref{eq:9IV}). 
C-R eq. (\ref{eq:7IV}) can also be cast in the form (\ref{eq:9IV}) or equivalently (\ref{eq:9III}-\ref{eq:10III}) with slight modification, namely

\begin{equation}\label{eq:10IV}
	\left( \partial_0, \vec{\nabla} \right) \circ \left( (N-1)U_0, \vec{U} \right) = 0.
\end{equation}

\vspace{0.4cm}
\textbf{2.} If for a linear transformation $L$ on $A$ there exists $a \in A = C, Q, \Theta$ such that

\begin{equation}\label{eq:11IV}
	\forall u \in A \ \ \ \ 
	L(u) = au
\end{equation}

\noindent then of course

\begin{equation}\label{eq:12IV}
	\forall u,q \ \ \ \
	\mathrm{Tr}\left\{ L(uq) \right\} = \mathrm{Tr} \left\{ (Lu)q \right\}
\end{equation}

Whether (\ref{eq:12IV}) implies for $L$ representability in the form (\ref{eq:11IV}) or not is an open question at present. A positive answer would provide us with an algebraic interpretation of the 
$A$-analyticity concept as introduced via (\ref{eq:4IV}) C-R eqs.

\vspace{0.4cm}
\textbf{3.} The formulation (\ref{eq:4IV}) of C-R eqs. (equivalent to ours for composition algebras) is appropriate for extension to any algebra with unit element though the question immediately arises, whether this extension is equivalent in some sense to the straightforward extension of (\ref{eq:10III}).

\vspace{0.2cm}
Finally, let us remark that the Fueter analyticity is a special case of Clifford analyticity . For a Clifford algebra the $\partial_{\overline{x}}$ operator may be regarded as a kind of square root of Laplace operator with respect to Clifford algebra multiplication. In this sense $\partial_\Theta$  and  $\partial_{\overline{\Theta}}$ operators are the "square roots" of the "$\circ$" product introduced in Section III. This product becomes a Clifford one for $A = \mathbf{Q}$.

\vspace{0.6cm}
\noindent \textbf{Acknowledgements}

This ten years later article was born in TeX due to kindness and abilities of my undergraduate
Student of computer science Maciej Dziemia\'nczuk. Thanks and glory Unto Him . 

The author thanks Prof. J. \L awrynowicz  for invitation to participate in "The Polish-Mexican Seminar on Generalized Cauchy -Riemann Structures and  Surface Properties of Crystals" - Kazimierz  Dolny ; August 98  , and for His kind hospitality there as well as for
discussions on the subject.

% % % % % % % Bibliography 

\end{document}